\documentclass[12pt]{article}
\usepackage{latexsym}
\usepackage{amssymb}
\usepackage[matrix,arrow]{xy}
\input{xy}
\xyoption{all}

\newtheorem{theorem}{THEOREM}[section]
\newtheorem{proposition}[theorem]{PROPOSITION}
\newtheorem{lemma}[theorem] {LEMMA}
\newtheorem{example}[theorem]{EXAMPLE}

\newtheorem{corollary}[theorem]{COROLLARY}

\newenvironment{proof}{\noindent {\it\bf Proof.} \rm}

\def\Hom{\mathop{\rm Hom}\nolimits}

\def\qed{\hfill \mbox{$\square$}}
\def\End{\mathop{\rm End}\nolimits}

\def\Ext{\mathop{\rm Ext}\nolimits}

\def\Ker{\mathop{\rm Ker}\nolimits}

\def\Im{\mathop{\rm Im}\nolimits}
\def\dim{\mathop{\rm dim_k}\nolimits}
\def\Der{\mathop{\rm Der}\nolimits}

\def\Int{\mathop{\rm Int}\nolimits}

\def\HH {\rm {HH}}
\def\SH {\rm {SH}}

\def\path{\mathop{\rightsquigarrow}\nolimits}
\def\mod{\mathop{\rm mod}\nolimits}

\def \Z{{\mathbb Z}}

\def\C{\mathcal C}
\def\D{\mathcal D}

\begin{document}

\title{The first Hochschild cohomology group of a schurian cluster-tilted algebra}

\author{Ibrahim Assem \footnote{D\'epartement de Math\'ematiques, Universit\'e de Sherbrooke, Sherbrooke, Qu\'ebec, Canada, J1K 2R1, E-mail address: ibrahim.assem@usherbrooke.ca} \ and Mar\'\i a Julia Redondo\footnote {Instituto de Matem\'atica, Universidad Nacional del Sur, Av. Alem 1253, (8000) Bah\'\i a Blanca, Argentina, E-mail address: mredondo@criba.edu.ar} \thanks{The second author is a researcher from CONICET, Argentina.} }

\date{}

\maketitle

\begin{abstract}
Given a cluster-tilted algebra $B$ we study its first Hochschild cohomology group ${\HH}^1(B)$ with coefficients in the $B$-$B$-bimodule $B$. We find several consequences when $B$ is representation-finite, and also in the case where $B$ is cluster-tilted of type $\tilde{\mathbb{A}}$.
\end{abstract}


\section{Introduction}
Cluster categories were introduced in \cite{BMRRT} and also in \cite{CCS} for type $\mathbb{A}$, in order to understand better the cluster algebras of Fomin and Zelevinsky \cite{FZ}.  Cluster-tilted algebras were defined in \cite{BMR} and also in \cite{CCS2} for type $\mathbb{A}$.  These algebras have been studied by several authors (see, for instance, \cite{ABS1,CCS2,BMR,BMR2}).  Our objective here is, for a cluster-tilted algebra $B$, to study its first Hochschild cohomology group ${\HH}^1(B)$ with coefficients in the $B$-$B$-bimodule $B$, see \cite{CE}. As a first step, we consider the case where $B$ is schurian: this includes the case of all representation-finite cluster-tilted algebras.  There are several reasons for this restriction.  Indeed, it was shown in \cite{ABS1} that, if $C$ is a tilted algebra, then the trivial extension $C \ltimes \Ext^2_C(DC,C)$ of $C$ by the $C$-$C$-bimodule $\Ext^2_C(DC,C)$ is cluster-tilted, and, conversely, any cluster-tilted algebra is of this form.  As a consequence, one can describe the ordinary quiver of $B$ knowing that of $C$.  But the relations are more difficult to compute.  However, they are known in the representation-finite case, see \cite{BMR2}, and one can assume, as an approximation, that the algebra is schurian.  Also, in a schurian algebra, normalized derivations are diagonalizable and this simplifies considerably our calculations.  The main step in our proof consists in defining an equivalence relation between the arrows in the quiver of $B$ which are not in the quiver of $C$.  The number of equivalence classes is then denoted as $n_{B,C}$.  Our main theorem is the following.

\begin{theorem}
Let $k$ be an algebraically closed field and $B$ be a schurian cluster-tilted algebra.  If $C$ is a tilted algebra such that $B=C \ltimes \Ext^2_C(DC,C)$ then there is a short exact sequence of abelian groups
$$0 \to {\HH}^1(C) \to {\HH}^1(B) \to k^{n_{B,C}} \to 0.$$
\end{theorem}

We find several interesting consequences in the representation-finite case. Also, if $B$ is a schurian cluster-tilted algebra, then we prove that ${\HH}^1(B)=0$ if and only if $B$ is hereditary with ordinary quiver a tree. We then also study the case where $B$ is cluster-tilted of type $\tilde{\mathbb{A}}$. In this case, $B$ is not necessarily schurian, but it is gentle (see \cite{ABCJP} or also \cite{BV}), hence we can apply the results of \cite{CS}.

We now describe the contents of the paper.  After a short preliminary section, we introduce our equivalence relation in Section 2, prove our main theorem and its consequences in Section 3, then study the case $\tilde \mathbb{A}$ in Section 4.

\section {Preliminaries.}
\subsection{Notations.}
Throughout this paper, algebras are basic, connected and finite dimensional over an algebraically closed field $k$.  For such an algebra $C$, there exists a quiver $Q$ and an admissible ideal $I$ of the path algebra $kQ$ of $Q$ such that $C \simeq kQ/I$.  We denote by $Q_0$ the set of points of $Q$, and by $Q_1$ its set of arrows. For a point $x \in Q_0$, we denote by $e_x$ the corresponding primitive idempotent of $C$.  For an arrow $\alpha \in Q_1$ we denote by $s(\alpha)$ and $t(\alpha)$, respectively, its source and its target, and by $\alpha^{-1}$ its formal inverse.  A \textit{walk} in $Q$ is a composition of arrows and formal inverses of arrows.  A \textit{relation} from $x \in Q_0$ to $y \in Q_0$ is a linear combination $\rho= \sum_{i=1}^{m} a_i w_i$ where each $a_i$ is a non-zero scalar and each $w_i$ is a path of length at least two from $x$ to $y$.  If $m=1$ then $\rho$ is \textit{monomial} and, if $m=2$ then it is \textit{binomial}.  Any admissible ideal of $kQ$ is generated by a finite set of relations.  Following \cite{BoG}, $C=kQ/I$ can equivalently be considered as a $k$-category with object class $C_0=Q_0$ and set of morphisms $C(x,y)$ from $x$ to $y$ equal to the quotient of the $k$-vector space $kQ(x,y)$ of all linear combinations of paths from $x$ to $y$ by the subspace $I(x,y) = I \cap kQ(x,y)$.  A full subcategory $C'$ of $C$ is \textit{convex} if any path in $C$ with source and target in $C'$ lies entirely in $C'$.  An algebra $C$ is \textit{triangular} if its quiver is acyclic, that is, without oriented cycles.  It is \textit{schurian} if, for any $x,y \in C_0$, we have $\dim C(x,y) \leq 1$.  If $C=kQ/I$ is schurian, then $I$ is generated by monomial and binomial relations.

By a $C$-module is meant a finitely generated right module.  We denote by $\mod C$ their category.  For more notations or facts about algebras or modules, we refer to \cite{ASS}.

\subsection{Cluster-tilted algebras.} Let $A$ be a hereditary algebra and $\D^b(\mod A)$ denote the bounded derived category over $\mod A$.  The \textit{cluster category} ${\C}_A$ is the orbit category $\D^b(\mod A)/F$ where $F=\tau^{-1}[1]$ is the composition of the Auslander-Reiten translation $\tau^{-1}$ with the shift $[1]$ of $\D^b(\mod A)$.  Then ${\C}_A$ is a triangulated category with almost split triangles.  An object $\tilde T$ is \textit{tilting} if $\Ext^1_{{\C}_A}(\tilde T, \tilde T)=0$ and the number of isomorphism classes of indecomposable summands of $\tilde T$ equals the rank of the Grothendieck group of $A$, see \cite{BMRRT}.  The endomorphism algebra of a tilting object $B=\End_{{\C}_A}(\tilde T)$ is a \textit{cluster-tilted algebra}, see \cite{BMR}.  It is shown in \cite{ABS1} that, if $T$ is a tilting module over a hereditary algebra $A$, so that $C=\End_A(T)$ is a tilted algebra, then the trivial extension $\tilde C= C \ltimes \Ext^2_C(DC,C)$ (the \textit{relation-extension} of $C$) is cluster-tilted and, conversely, any cluster-tilted algebra is of this form (but in general, not uniquely: see \cite{ABS2}).  As a consequence, we have a description of the quiver of $\tilde C$.  Let $R$ be a \textit{system of relations} for the tilted algebra $C=kQ/I$, that is, $R$ is a subset of $\cup_{x,y \in Q_0} I(x,y)$ such that $R$, but no proper subset of $R$, generates $I$ as an ideal of $kQ$.  It is shown in \cite{ABS1} that the quiver  $\tilde Q$ of $\tilde C$ is as follows:
\begin{itemize}
\item[(a)] $\tilde Q_0 = Q_0$;
\item[(b)] For $x,y \in Q_0$, the set of arrows in $\tilde Q$ from $x$ to $y$ equals the set of arrows in $Q$ from $x$ to $y$ (which we call \textit{old arrows}) plus $\vert R \cap I(y,x)\vert$ additional arrows (which we call \textit{new arrows}).
\end{itemize}

\begin{lemma}\label{13}
Let $C=kQ/I$ be a schurian tilted algebra, then its relation-extension $\tilde C$ contains no walk of the form $w=\alpha w' \beta$, where $\alpha, \beta$ are new arrows, and $w'$ is a walk not containing zero relations and
consisting entirely of old arrows.
\end{lemma}

\begin{proof}
Suppose there exists such a walk, and assume, without loss of generality, that the length of $w'$ is minimal.  Since new arrows correspond to relations in $C$, and the quiver $Q$ is acyclic, then the existence of such a walk in the quiver $\tilde Q$ of $\tilde C$ implies that $C$ contains a subquiver (maybe not full) of one of the forms
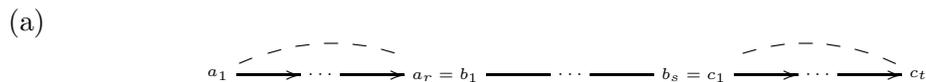
\begin{itemize}
\item [(a)]
\begin{tiny}
\[ \xymatrix{ a_1  \ar[r] \ar@/^1pc/@{--}[rr]     & \cdots  \ar[r] & a_r=b_1 \ar@{-}[r] & \cdots \ar@{-}[r] & b_s=c_1 \ar[r] \ar@/^1pc/@{--}[rr] & \cdots \ar[r] &  c_t
} \]
\end{tiny}

\item [(b)] 
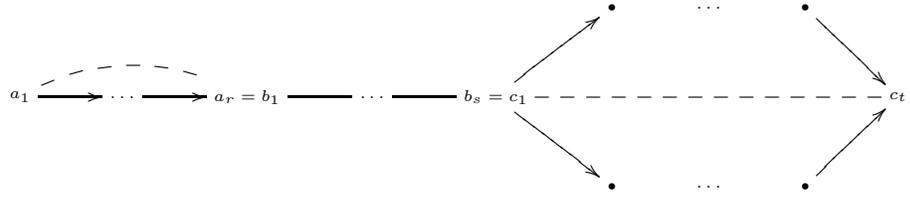
\begin{tiny}
\[ \xymatrix{  &&&&& \scriptscriptstyle \bullet & \dots & \scriptscriptstyle \bullet \ar[dr] \\
a_1  \ar[r] \ar@/^1pc/@{--}[rr]    & \dots \ar[r] & a_r=b_1 \ar@{-}[r] & \dots \ar@{-}[r] & b_s=c_1 \ar[ur] \ar[dr] \ar@{--}[rrrr] & & &  & c_t \\
&&&&& \scriptscriptstyle \bullet & \dots & \scriptscriptstyle \bullet \ar[ur] } \]
\end{tiny} 

\item [(c)]
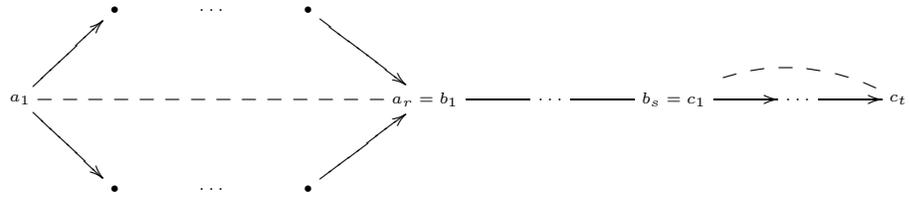
\begin{tiny}
\[ \xymatrix{  & \scriptscriptstyle \bullet & \dots & \scriptscriptstyle \bullet \ar[dr]   \\
a_1  \ar[ur] \ar[dr] \ar@{--}[rrrr] & & &  & a_r=b_1 \ar@{-}[r]    & \dots \ar@{-}[r] & b_s=c_1 \ar[r] \ar@/^1pc/@{--}[rr]    & \dots \ar[r] & c_t \\
& \scriptscriptstyle \bullet & \dots & \scriptscriptstyle \bullet \ar[ur]  } \]
\end{tiny}

\item [(d)]
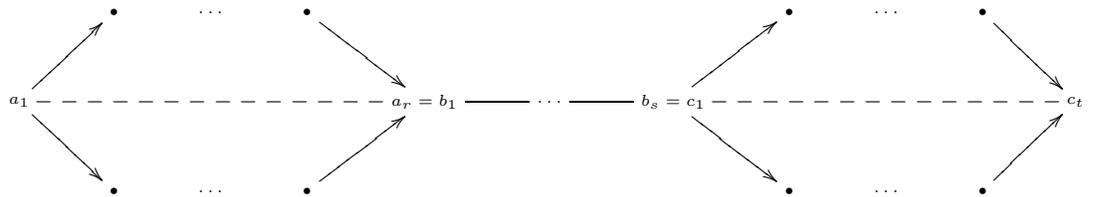
\begin{tiny}
\[ \xymatrix{  & \scriptscriptstyle \bullet &  \dots & \scriptscriptstyle \bullet \ar[dr] &&& & \scriptscriptstyle \bullet & \dots & \scriptscriptstyle \bullet \ar[dr] \\
a_1  \ar[ur] \ar[dr] \ar@{--}[rrrr] & & &  & a_r=b_1 \ar@{-}[r]    & \dots \ar@{-}[r] & b_s=c_1 \ar[ur] \ar[dr] \ar@{--}[rrrr] & & &  & c_t \\
& \scriptscriptstyle \bullet & \dots & \scriptscriptstyle \bullet \ar[ur] &&& & \scriptscriptstyle \bullet & \dots & \scriptscriptstyle \bullet \ar[ur] } \]
\end{tiny}

\item [(e)]
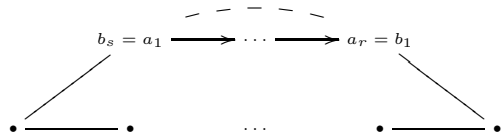
\begin{tiny}
\[ \xymatrix{ & b_s=a_1 \ar@/^1pc/@{--}[rr] \ar[r] & \dots \ar[r] & a_r=b_1 \ar@{-}[dr] \\
\scriptscriptstyle \bullet \ar@{-}[r] \ar@{-}[ru]  & \scriptscriptstyle \bullet & \dots & \scriptscriptstyle \bullet \ar@{-}[r] & \scriptscriptstyle \bullet} \]
\end{tiny}

\item [(f)]
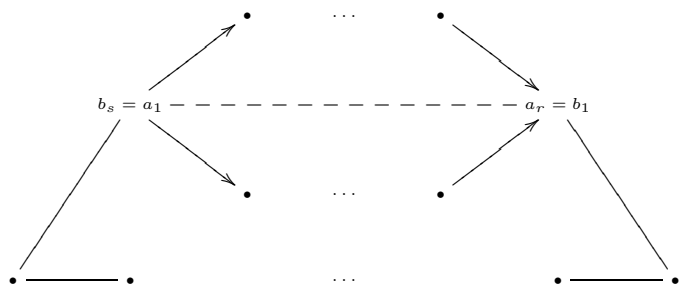
\begin{tiny}
\[ \xymatrix{  & & \scriptscriptstyle \bullet & \dots & \scriptscriptstyle \bullet \ar[dr] \\
& b_s=a_1 \ar[ur] \ar[dr] \ar@{--}[rrrr] & & &  & a_r=b_1 \ar@{-}[ddr] \\
& & \scriptscriptstyle \bullet & \dots & \scriptscriptstyle \bullet \ar[ur] \\
\scriptscriptstyle \bullet \ar@{-}[r] \ar@{-}[ruu]  & \scriptscriptstyle \bullet & & \dots &  & \scriptscriptstyle \bullet \ar@{-}[r] & \scriptscriptstyle \bullet} \]
\end{tiny}
\end{itemize}
where the walk $b_1 -  \cdots - b_s$ (with $s \geq 1$) is non-zero and dotted lines represent relations.  Here the last two cases occur in case $\alpha=\beta$.  Let $C'$ be the full subcategory of $C$ generated by the points $a_i, b_j, c_k$. By \cite[III.6.5, p. 146]{Ha1}, $C'$ is a tilted algebra. It follows from our minimality assumption and the hypothesis that $\tilde C$, and hence $C$, are schurian, that in $C'$ there is no additional arrow between two ${b_j}$'s.  Now, in each of the cases above, let $M$ be the $C'$-module defined as a representation by
\[ M(x)= \left\{ \begin{array}{lll}
k, &\mbox{if $x \in \{ b_1, \dots, b_s\}$}, \\
0, & \mbox{if $x \not \in \{ b_1, \dots, b_s\}$,}
\end{array} \right. \]
also, if $\alpha$ is an arrow such that $s(\alpha)$ and $t(\alpha)$ both belong to $\{ b_1, \dots, b_s\}$, we let $M(\alpha)=1$, while, for all other arrows $\alpha$, we let $M(\alpha)=0$.  Then $M$ is an indecomposable $C'$-module, and it is can be seen that both its projective and its injective dimensions equal $2$, a contradiction because $C'$ is tilted. \qed
\end{proof}

\paragraph{} In this paper we mainly consider schurian cluster-tilted algebras.  We now prove that re\-pre\-sen\-tation-finite cluster-tilted algebras are schurian.  We recall that, if $B$ is a representation-finite cluster-tilted algebra and $C$ is a tilted algebra such that $B=\tilde C$, then $C$ is tilted of Dynkin type, see \cite{BMR}.

\begin{lemma} Let $B$ be a representation-finite cluster-tilted algebra.  Then $B$ is schurian.
\end{lemma}

\begin{proof} 
Let $x,y \in B_0$ be such that $\dim B(x,y) \geq 2$ and $C$ be a tilted algebra such that $B = \tilde C$.  Assume first that $x=y$.  Then $\dim B(x,x) \geq 2$.  In particular, there exists in the quiver $\tilde Q$ of $B$ a non-zero cycle $\gamma$ from $x$ to $x$.  This cycle $\gamma$ must contain a new arrow, because $C$ is triangular.  But then, by Lemma \ref{13}, such an arrow must be unique, and this yields a contradiction, again by Lemma \ref{13}.  Therefore $x \not=y$.  Let $e=e_x+e_y$, then $eBe$ has as quiver

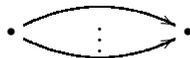
\[ \xymatrix{\scriptscriptstyle \bullet \ar@/^.8pc/[rr] \ar@/_.8pc/[rr] & \vdots & \scriptscriptstyle \bullet   } \]
and thus is representation-infinite, a contradiction to the representation-finiteness of $B$. \qed
\end{proof}

\subsection {Example} Clearly, there exist schurian representation-infinite cluster-tilted algebras.  Let $C$ be given by the quiver
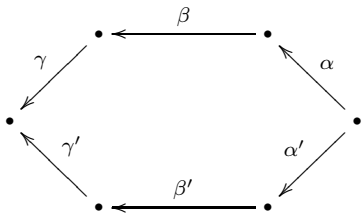
\[ \xymatrix{ & \ar[dl]_\gamma \scriptscriptstyle \bullet &  & \ar[ll]_\beta \scriptscriptstyle \bullet  \\
\scriptscriptstyle \bullet & & & & \scriptscriptstyle \bullet \ar[dl]_{\alpha'} \ar[ul]_\alpha\\
& \ar[ul]_{\gamma'} \scriptscriptstyle \bullet & & \scriptscriptstyle \bullet  \ar[ll]_{\beta'}   } \]
bound by $\alpha \beta =0, \alpha'\beta'=0$. Then $C$ is tilted of type $\tilde \mathbb{A}$ and $B=\tilde C$ is given by the quiver
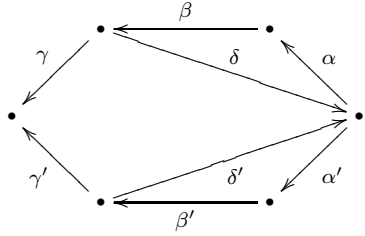
\[ \xymatrix{ & \ar[dl]_\gamma \ar[drrr]^\delta \scriptscriptstyle \bullet &  & \ar[ll]_\beta \scriptscriptstyle \bullet  \\
\scriptscriptstyle \bullet & & & & \scriptscriptstyle \bullet \ar[dl]^{\alpha'} \ar[ul]_\alpha\\
& \ar[ul]^{\gamma'} \ar[urrr]_{\delta'} \scriptscriptstyle \bullet & & \scriptscriptstyle \bullet  \ar[ll]^{\beta'}   }\]
bound by $\alpha \beta =0, \beta\delta=0, \delta \alpha=0,  \alpha'\beta'=0,\beta'\delta'=0, \delta'\alpha'=0$, see \cite{ABCJP} or Section \ref{41} below.  Thus, $B$ is representation-infinite and schurian.

\section{Arrow equivalence}

\begin{lemma}\label{21}
Let $B=k \tilde Q / \tilde I$ be a schurian cluster-tilted algebra, and $C=kQ/I$ be a tilted algebra such that $B=\tilde C$.  Let $\rho=a_1 w_1+a_2w_2$ be a binomial relation in $\tilde I$.  Then either $\rho$ is a relation in $I$, or there exist exactly two new arrows $\alpha_1, \alpha_2$ such that $w_1=u_1\alpha_1v_1, w_2=u_2\alpha_2v_2$ (with $u_1, u_2, v_1, v_2$ paths consisting entirely of old arrows).
\end{lemma}

\begin{proof}
By the main result of \cite[2.4]{ACT} (see also \cite[3.2]{ACMT}), either $\rho$ is a relation in $I$, or there exist at least two new arrows $\alpha_1, \alpha_2$ such that $w_1=u_1\alpha_1v_1, w_2=u_2\alpha_2v_2$, with $u_1, u_2, v_1, v_2$ paths.  The uniqueness of each of $\alpha_1, \alpha_2$ follows from Lemma \ref{13}. \qed
\end{proof}

\medskip

The Lemma \ref{21} above brings us to our main definition.  Let $B=k \tilde Q /\tilde I$ be a schurian cluster-tilted algebra, and $C=kQ/I$ be a tilted algebra such that $B=\tilde C$.  We define a relation $\sim$ on the set $\tilde Q_1 \setminus Q_1$ of new arrows as follows. For every $\alpha \in \tilde Q_1 \setminus Q_1$, we set $\alpha \sim \alpha$.  If $\rho=a_1 w_1+ a_2 w_2$ is a binomial relation in $\tilde I$ and $\alpha_1, \alpha_2$ are as in Lemma \ref{21} above, then we set $\alpha_1 \sim \alpha_2$.

By Lemma \ref{13}, the relation $\sim$ is unambiguously defined.  It is clearly reflexive and symmetric.  We let $\approx$ be the least equivalence relation defined on the set $\tilde Q_1 \setminus Q_1$ such that $\alpha \sim \beta$ implies $\alpha \approx \beta$ (that is, $\approx$ is the transitive closure of $\sim$).

We define the \textit{relation invariant} of $B$ to be the number $n_{B,C}$ of equivalence classes under the relation $\approx$.

We now prove that, in the representation-finite case, the two relations $\sim$ and $\approx$ coincide (so, $\sim$ is an equivalence relation).

\begin{lemma}
Let $B= k \tilde Q/\tilde I$ be a representation-finite schurian cluster-tilted algebra, and $C=kQ/I$ be a tilted algebra such that $B=\tilde C$.  Then there exists no triple $(\alpha, \beta, \gamma)$ of pairwise distinct new arrows such that $\alpha \sim \beta$ and $\beta \sim \gamma$.
\end{lemma}

\begin{proof}
Since every new arrow of $\tilde Q$ corresponds to a relation in $I$, then there exist three relations $\rho, \sigma, \xi$ in $I$ such that each of the pairs $(\rho, \sigma)$ and $(\sigma, \xi)$ belongs to two parallel paths.

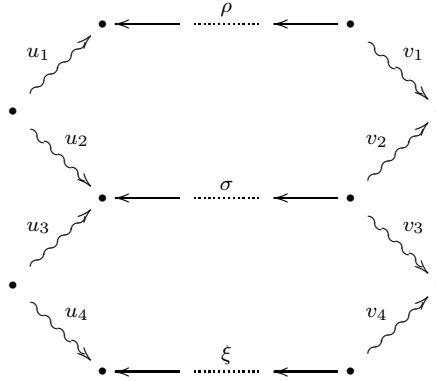
\[ \xymatrix{ & \scriptscriptstyle \bullet & \ar[l]  \ar@{.}[r]^\rho &  & \ar[l] \scriptscriptstyle \bullet \ar@{~>}[dr]^{v_1} & \\
\scriptscriptstyle \bullet \ar@{~>}[ur]^{u_1} \ar@{~>}[dr]^{u_2} & & &  & & \scriptscriptstyle \bullet \\
& \scriptscriptstyle \bullet & \ar[l] \ar@{.}[r]^\sigma &  & \ar[l]  \scriptscriptstyle \bullet \ar@{~>}[ur]^{v_2} \ar@{~>}[dr]^{v_3}& \\
\scriptscriptstyle \bullet \ar@{~>}[ur]^{u_3} \ar@{~>}[dr]^{u_4} & & & & & \scriptscriptstyle \bullet \\
& \scriptscriptstyle \bullet & \ar[l] \ar@{.}[r]^\xi  &  & \ar[l]  \scriptscriptstyle \bullet \ar@{~>}[ur]^{v_4} & }\]
By Lemma \ref{13}, we get that $u_i, v_j$ are actually paths in $C$.  If at least two of the relations $\rho, \sigma, \xi$ are binomial, then the algebra $C$ contains a representation-infinite full subcategory.  If exactly one of these three relations is binomial, then we obtain a contradiction to Lemma \ref{13} again.  Thus, we may assume the three relations to be monomial.

Let $E$ be the convex envelope in $C$ of all the points belonging to one of the paths $u_i, v_j, \rho, \sigma, \xi$.  We claim that $E$ does not contain a non-trivial path from a point in $u_i$ to a point in $u_j$, with $i \not = j$. Indeed, if such a path exists and contains a relation, then we get a contradiction to Lemma \ref{13}.  If it contains no relation and is of length at least one, then we have a representation-infinite full subcategory.  This establishes our claim.\\
On the other hand, there may exist a non-trivial path from a point in $u_i$ to a point in $v_j$.  In this case, this path must contain a monomial relation, because of Lemma \ref{13}.  Let $d$ be the total number of relations between the points of $E$, then $d \geq 3$.  Let $v=(v_x)_{x \in E_0}$ be the following vector
\[ v_x= \left\{ \begin{array}{rll}
-1, &\mbox{if $x$ lies in the convex envelope of the $u_i$}, \\
1, &\mbox{if $x$ lies in the convex envelope of the $v_j$}, \\
0, & \mbox{otherwise.}
\end{array} \right. \]
Evaluating the Tits quadratic form of $E$ in this vector yields 
\[ \begin{array}{ll}
q_E(v) &= \sum v_x^2 - \sum_{x \stackrel{\alpha}{\to} y} v_x v_y + \sum_{x \stackrel{\rho}{\to} y} v_x v_y\\
&= 1 + 1 - d \\
&= 2-d <0
\end{array} \]
In particular, it is not positive definite.  This contradicts the fact that $E$ is a full convex subcategory of $C$, which is tilted of Dynkin type. \qed
\end{proof}

\begin{corollary}
Let $B=k \tilde Q / \tilde I$ be a representation-finite cluster-tilted algebra, and $C=kQ/I$ be a tilted algebra such that $B= \tilde C$.  Let $R$ and $\tilde R$ be respectively systems of relations for $I$ and $\tilde I$, and $n'_{B,C}$ be the number of binomial relations which belong to $\tilde R \setminus R$, then $n_{B,C} = |R|- n'_{B,C}$.
\end{corollary}

\section{The main result and its consequences}

For the Hochschild cohomology groups, we refer the reader to \cite{CE,Ha2}.  We recall that, it $B$ is a finite dimensional $k$-algebra, then the first Hochschild cohomology group of $B$ can be written as ${\HH}^1(B)= \Der_0 B / \Int_0 B$, where $\Der_0 B$ is the $k$-vector space of the \textit{normalized derivations}, that is, of the derivations $\delta$ such that $\delta(e_x)=0$ for every primitive idempotent $e_x$, and $\Int_0 B$ is the subspace of the interior normalized derivations (see, for instance, \cite[3.1]{Ha2}).

\begin{theorem}\label{31}
Let $B$ be a schurian cluster-tilted algebra, and $C$ be a tilted algebra such that $B=\tilde C$.  Then there exists a short exact sequence of abelian groups
$$0 \to {\HH}^1(C) \to {\HH}^1(B) \to k^{n_{B,C}} \to 0.$$
\end{theorem}

\begin{proof}
As usual, we let $B=k\tilde Q / \tilde I$ and $C=kQ/I$.  Let $\delta \in \Der_0 B$ then, for every arrow $\alpha: x \to y$ in $\tilde Q_1$, we have $\alpha=e_x \alpha e_y$, hence 
$$\delta(\alpha)=\delta(e_x \alpha e_y)= e_x \delta(\alpha) e_y \in e_x B e_y.$$
Since $B$ is schurian, there exists a scalar $\lambda_\alpha$ such that $\delta (\alpha)= \lambda_\alpha \alpha$. Let then $\alpha_1 \alpha_2 \dots \alpha_t$ be a path in $\tilde Q$, then
$$\delta (\alpha_1 \alpha_2 \dots \alpha_t)= (\lambda_{\alpha_1} + \dots + \lambda_{\alpha_t}) \alpha_1 \alpha_2 \dots \alpha_t.$$
In particular, $\delta$ is uniquely determined by its value on the arrows. \\
Let $\rho= a \ \beta_1 \dots \beta_r - b \ \gamma_1 \dots \gamma_s$, with $a,b \in k^*$, be a binomial relation, then $\rho \in \tilde I$ implies 
$$0= a \ \delta (\beta_1 \dots \beta_r) - b \ \delta (\gamma_1 \dots \gamma_s)$$
which yields
$$a \ (\sum_{i=1}^r \lambda_{\beta_i}) \ \beta_1 \cdots \beta_r = b \ (\sum_{j=1}^s \lambda_{\gamma_j})\ \gamma_1 \cdots \gamma_s$$
and hence 
$$\sum_{i=1}^r \lambda_{\beta_i} = \sum_{j=1}^s \lambda_{\gamma_j}.$$
We use this observation to define a morphism $\zeta: \Der_0 C \to \Der_0 B$ as follows.  Let $\delta \in \Der_0 C$, then we let $\tilde \delta= \zeta (\delta)$ be defined by its action on the arrows according to 
\[\tilde \delta (\alpha) = \left\{ 
\begin{array}{rll}
\delta(\alpha),  &\mbox{if $\alpha \in Q_1$},\\
- (\sum_{i=1}^r \lambda_{\beta_i}) \alpha, & \mbox{if $\alpha \in \tilde Q_1 \setminus Q_1$},
\end{array} \right. \]
where $\beta_1 \cdots \beta_r$ is a path appearing in the relation in $I$ which defines the new arrow $\alpha$.  Since $B=\tilde C$ is schurian, this relation is either monomial or binomial and, in the latter case, the above argument shows that $\sum_{i=1}^r \lambda_{\beta_i}$ is uniquely determined and hence $\tilde \delta (\alpha)$ is unambiguously defined.  Clearly, $\zeta: \Der_0 C \to \Der_0 B$ is an injective map.\\
We now claim that $\zeta (\Int_0 C)=\Int_0 B$.  Let $\delta_a \in \Int_0 C$, then there exists $a \in C$ such that $\delta_a (c)=ac-ca$ (for every $c \in C$).  Since $\delta_a$ is normalized, we have, for every primitive idempotent $e_x$,
$$0=\delta_a(e_x)=ae_x-e_xa.$$
Hence 
\[ a=a.1=a \sum e_x = \sum a e_x = \sum (a e_x) e_x = \sum (e_x a) e_x = \sum a_x e_x \]
with $a_x \in k$ for every $x$.  This indeed follows from the fact that $e_x a e_x \in e_x C e_x$ and $C$ is schurian.\\
Let $\tilde \delta = \zeta (\delta_a)$.  For $\alpha \in Q_1$, we have
$$\tilde \delta (\alpha) = \delta_a(\alpha) = a \alpha - \alpha a = (a_{s(\alpha)}-a_{t(\alpha)}) \alpha.$$
Let now $\alpha \in \tilde Q_1 \setminus Q_1$ and $\beta_1 \cdots \beta_r$ be a path appearing in the relation defining $\alpha$.  Note that
$$\delta_\alpha (\beta_i)= \lambda_{\beta_i} \beta_i = (a_{s(\beta_i)}-a_{t(\beta_i)}) \beta_i$$
for every $i$ such that $1 \leq i \leq r$.  Therefore
$$\sum_{i=1}^r \lambda_{\beta_i} = \sum_{i=1}^r (a_{s(\beta_i)}-a_{t(\beta_i)}) = - a_{s(\alpha)} + a_{t(\alpha)}$$
so we have
$$ \tilde \delta (\alpha) = - (\sum_{i=1}^r \lambda_{\beta_i}) \alpha = (a_{s(\alpha)} - a_{t(\alpha)}) \alpha.$$
This shows that $\zeta(\Int_0 C) \subseteq \Int_0 B$.  Now, because $\delta_a$ is determined by the element $a=\sum a_x e_x$ and $Q_0=\tilde Q_0$, we actually have $\zeta(\Int_0 C)=\Int_0 B$.  This establishes our claim.\\
We now define a map
$$\phi: \Der_0 B \to k^{n_{B,C}}$$
as follows.  For a derivation $\delta \in \Der_0 B$, let $\delta |_C$ denote its restriction to $C$.  Clearly, $\delta |_C \in \Der_0 C$.  We let $\widetilde{\delta |_C}= \zeta( \delta |_C)$ then we set
$$\phi (\delta) = (\lambda_\alpha - \sum_{i=1}^r \lambda_{\beta_i})_{\alpha \in S}$$
where $S$ is a complete set of representatives of the equivalence classes of the new arrows under the relation $\approx$, and, for $\alpha \in S$, we let $\beta_1 \cdots \beta_r$ be a path in $C$ occurring in the relation defining $\alpha$. As observed above, the sum $\sum_{i=1}^r \lambda_{\beta_i}$ is uniquely determined by $\alpha$.  We still have to prove that $\phi (\delta)$ does not depend on the particular representative chosen in the class of $\alpha$.  Assume then that $\alpha \sim \alpha'$.  Then there exists a binomial relation $\rho= a(u\alpha v)- a' (u' \alpha' v')$, with $a,a' \in k^*$ and $u,v,u',v'$ paths in $C$.  Let $\beta'_1 \cdots \beta'_s$ be a path occurring in a relation defining $\alpha'$.  Then we have 
\[ \begin{array}{lll}
(\lambda_\alpha - \sum_{i=1}^r \lambda_{\beta_i}) a (u \alpha v) & = & (\delta - \widetilde{\delta|_C})(a u \alpha v)\\
& = & (\delta - \widetilde{\delta|_C})(a' u' \alpha' v')\\
& = & (\lambda_{\alpha'} - \sum_{j=1}^s \lambda_{\beta'_j}) a' (u' \alpha' v').
\end{array}\]
This shows that $\phi$ is unambiguously defined.  Clearly, $\phi$ is $k$-linear.\\
We now show that $\phi$ is surjective.  Let $(\mu_\alpha)_{\alpha \in S} \in k^{n_{B,C}}$ and define $\delta \in \Der_0 B$ by its value on the arrows as follows
\[ \delta (\alpha) = \left \{ \begin{array}{lll}
0 & \mbox{if $\alpha \in Q_1$}\\
\mu_{\alpha'} \alpha & \mbox{ if $\alpha \in \tilde Q_1 \setminus Q_1$ and $\alpha \sim \alpha'$}
\end{array} \right. \]
This is clearly a derivation and moreover $\delta|_C=0$.  Then $\phi(\delta)=(\mu_\alpha)_{\alpha \in S}$.\\
Finally, we prove that $\Ker \phi = \Der_0 C$.  Indeed, $\Ker \phi$ consists of the $\delta \in \Der_0 B$ such that $\delta - \widetilde{\delta|_C}=0$.  By definition, this set equals $\Im \zeta = \Der_0 C$.\\
We have thus shown that there exists a short exact sequence of $k$-vector spaces
$$0 \to \Der_0 C \stackrel{\zeta}{\to} \Der_0 B \stackrel{\phi}{\to} k^{n_{B,C}} \to 0.$$
Since $\zeta(\Int_0 C)= \Int_0 B$, the statement follows at once. \qed
\end{proof}

\medskip

Before our first corollary, we introduce some notation.  Let $B=k\tilde Q/\tilde I$ be a schurian cluster-tilted algebra.  By \cite{BM}, the fundamental group $\pi_1(\tilde Q, \tilde I)$ does not depend on the particular presentation of $B$, we may then denote it as $\pi_1(B)$.  If $B$ is representation-finite, then $\pi_1(B)$ is a free group \cite{MP}.  We denote by $L_m$ the free group on $m$ letters.  Also, let ${\SH}_1(B)$ denote the first (simplicial) homology group of the classifying space of $B$, see \cite{Bu,BrG}.

\begin{corollary}\label{32}
Let $B$ be a representation-finite cluster-tilted algebra.  Then: 
\begin{itemize}
\item[(a)] ${\HH}^1(B)=k^{n_{B,C}}$;
\item[(b)] $\pi_1(B)\cong L_{n_{B,C}}$;
\item[(c)] ${\SH}^1(B) \cong {\Z}^{n_{B,C}}$;
\item[(d)] $|R| \geq n'_{B,C}$ for every tilted algebra $C$ such that $B=\tilde C$.
\end{itemize}
\end{corollary}

\begin{proof}
\begin{itemize}
\item []
\item[(a)] Since $C$ is tilted of Dynkin type, then ${\HH}^1(C)=0$.  The result follows then from Theorem \ref{31}.
\item[(b)] Since $B$ is schurian, we have ${\HH}^1(B) \cong \Hom(\pi_1(B), k^+)$ by \cite{SP}.  Since $\pi_1(B) \cong L_m$ for some $m$, comparing dimensions yields $m=n_{B,C}$.
\item[(c)] This follows from the Poincar\'e-Hurwitz theorem.
\item[(d)] This is trivial. \qed
\end{itemize}
\end{proof}

\begin{corollary}
Let $B$ be a schurian cluster-tilted algebra.  Then ${\HH}^1(B)=0$ if and only if $B$ is hereditary having a tree as ordinary quiver.
\end{corollary}

\begin{proof}
Since the sufficiency is obvious, we prove the necessity. Let $B$ be a schurian cluster-tilted algebra and $C$ be a tilted algebra such that $B=\tilde C$. By Theorem \ref{31}, ${\HH}^1(B)=0$ implies $n_{B,C}=0$.  But then there are no new arrows. Consequently, $B=C$ and it is therefore hereditary.  By \cite{Ha2} the quiver of $B$ is a tree.
\qed
\end{proof}

\begin{corollary}\label{4.3}
Assume $B$ is a schurian cluster-tilted monomial algebra.  Then ${\HH}^1(B) \cong {\HH}^1(C) \oplus k^{|R|}$, as $k$-vector spaces, for every tilted algebra $C$ such that $B=\tilde C$.  If $C$ is monomial and ${\HH}^1(B) \cong {\HH}^1(C) \oplus k^{|R|}$, then $B$ is monomial.
\end{corollary}

\begin{proof}
The first statement follows directly from Theorem \ref{31}.  Assume now that $C$ is monomial but $B$ is not.  Then there exists a binomial relation in $B$ but none in $C$. Then $n'_{B,C}\geq 1$ and ${\HH}^1(B) \cong {\HH}^1(C) \oplus k^{n_{B,C}}$ gives a contradiction by comparing dimensions. \qed
\end{proof}

\subsection{Example}
The dimension of ${\HH}^1(B)$ is unbounded, even in the representation-finite case.  Let $C_1$ be the tilted algebra given by the quiver

\[ \xymatrix{ \scriptscriptstyle \bullet \ar[r]^{\alpha_1} & \scriptscriptstyle \bullet \ar[r]^{\beta_1} & \scriptscriptstyle \bullet } \]
bound by $\alpha_1 \beta_1 =0$ and, for every $d \geq 1$, let $C_d$ be given by the quiver

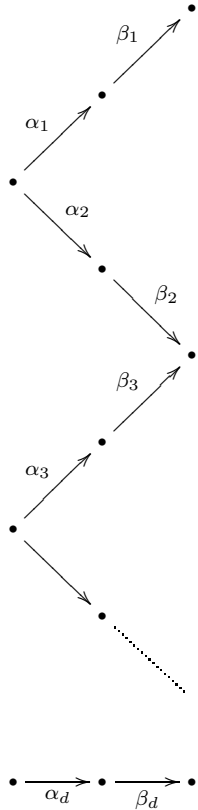
\[ \xymatrix{ & & \scriptscriptstyle \bullet \\
& \scriptscriptstyle \bullet \ar[ur]^{\beta_1}& \\
\scriptscriptstyle \bullet \ar[ur]^{\alpha_1} \ar[dr]^{\alpha_2} & & \\
& \scriptscriptstyle \bullet \ar[dr]^{\beta_2} & \\
& & \scriptscriptstyle \bullet \\
& \scriptscriptstyle \bullet \ar[ur]^{\beta_3} & \\
\scriptscriptstyle \bullet \ar[ur]^{\alpha_3} \ar[dr] & & \\
& \scriptscriptstyle \bullet \ar@{.}[dr] & \\ 
& & \\
\scriptscriptstyle \bullet \ar[r]_{\alpha_d} & \scriptscriptstyle \bullet \ar[r]_{\beta_d} & \scriptscriptstyle \bullet } \]
bound by $\alpha_i \beta_i =0$ for every $i$.  Then $C_d$ is tilted of Dynkin type $\mathbb{A}$ and by \cite{ABCJP} or \cite{BV}, $B$ is monomial.  Applying Corollary \ref{4.3}, we get ${\HH}^1(B)=k^d$.

\section{Cluster-tilted algebras of type $\tilde \mathbb{A}$} \label{41}

If $B$ is a cluster-tilted of type $\tilde{\mathbb{A}}$ then it is gentle, and hence monomial, see \cite{ABCJP}. 
We recall from \cite{ASk} that an algebra $C=kQ/I$ is called \textit{gentle} if:
\begin{itemize}
\item[(1)] For any $x \in Q_0$, there are at most two arrows of $Q$ having $x$ as a source or as a target;
\item[(2)] $I$ is generated by paths of length two;
\item[(3)] For any $\alpha \in Q_1$, there is at most one arrow $\beta$ and one arrow $\gamma$ such that $\alpha \beta \not \in I$ and $\gamma \alpha \not \in I$;
\item[(4)] For any $\alpha \in Q_1$, there is at most one arrow $\xi$ and one arrow $\eta$ such that $\alpha \xi \in I$ and $\eta \alpha \in I$.
\end{itemize}
As observed in \cite{ABCJP}, if $C$ is tilted and gentle, then either it is of type $\mathbb{A}$ (in which case $Q$ is a gentle tree containing no double zeros, that is, walks of the form $w=\alpha \beta w' \gamma \delta$ where $\alpha, \beta, \gamma, \delta \in Q_1$ are such that $\alpha \beta, \gamma \delta \in I$, and $w'$ is a non-zero walk), or else of type $\tilde \mathbb{A}$ (in which case it contains a unique non-oriented cycle without double zeros and all arrows attached to the cycle either enter it or all leave it).  It is shown there that $B=\tilde C$ is also gentle and in fact all relations occur in $3$-cycles, that is, cycles of the form

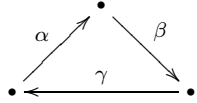
\[ \xymatrix{ & \ar[dr]^\beta \scriptscriptstyle \bullet \\
\ar[ur]^\alpha \scriptscriptstyle \bullet & & \ar[ll]_\gamma \scriptscriptstyle \bullet  }\] 
where $\alpha \beta, \beta \gamma$ and $\gamma \alpha$ are relations in $\tilde C$.

\begin{proposition}
Let $B$ be a cluster-tilted algebra of type $\tilde \mathbb{A}$, and $C$ be a tilted algebra such that $B=\tilde C$.  Let $R$ be a system of relations for $C$, then $${\HH}^1(B) = k^{|R|+\epsilon +1}$$ where $0 \leq \epsilon \leq 2$.
\end{proposition}

\begin{proof}
Since $B$ is monomial, it follows from \cite{CS} that, in their notation,
$$\dim {\HH}^1(B)= \dim Z(B)- |\tilde Q_0 ||N| + | \tilde Q_1 ||N| - |(\tilde Q_1 || N)_e| - \dim \Im R_g.$$
We now study each of these terms.  Firstly, $Z(B)$ is the center of $B$, so $\dim Z(B)=1$.  Secondly, $\tilde Q_0 || N$ is the set of non-zero oriented cycles in $(\tilde Q, \tilde I)$ (where, as usual, $B=k \tilde Q/\tilde I$), including the points.  Then
$$|\tilde Q_0 ||N| = |\tilde Q_0|=|Q_0|.$$
Thirdly, $\tilde Q_1 ||N$ is the set of pairs consisting of an arrow, and a non-zero path (including points) parallel to this arrow. Since $Q$ contains a unique non-oriented cycle, we have 
$$|\tilde Q_1 ||N|= |\tilde Q_1| + \epsilon'$$
where $\epsilon'=2$ whenever $\tilde Q$ contains a double arrow, $\epsilon'=1$ whenever $\tilde Q$ contains a bypass (that is, an arrow $\alpha$ such that there exists a path $\beta_1 \cdots \beta_s$ parallel to $\alpha$) and $\epsilon'=0$, otherwise. Next,
$$(\tilde Q_1 || N)_e = (\tilde Q_1 || N) \setminus \{ (\tilde Q_1 || N)_g \cup (\tilde Q_1 || N)_a \}$$
where
\begin{itemize}
\item[(1)] $(\tilde Q_1 || N)_g$ is the set of pairs $(\alpha, \gamma) \in \tilde Q_1 || N$ where $\gamma$ is either a point, or a path starting or ending with the arrow $\alpha$.  In view of the discussion above, the only such pairs are of the form $(\alpha, \alpha)$ with $\alpha \in \tilde Q_1$.  Therefore $|(\tilde Q_1 || N)_g| = | \tilde Q_1|$.
\item[(2)] $(\tilde Q_1 || N)_a$ is the set of pairs $(\alpha, \gamma) \in \tilde Q_1 || N$ such that, in each relation $\rho$ where $\alpha$ appears, replacing $\alpha$ by $\gamma$ yields a zero path. 
\end{itemize}
Finally, let $R||N$ be the set of pairs $(w,w')$ where $w,w'$ are parallel paths such that $w \in R$ and $w'$ is non-zero.  Then $R_g$ is the linear map from $k(\tilde Q_1 || N)_g$ to $k(R||N)$ defined as follows.  Let $\rho$ be a relation in which a given arrow $\alpha$ appears and assume $(\alpha, \gamma) \in (\tilde Q_1 || N)_g$, then denote by $\rho'$ a non-zero parallel path obtained from $\rho$ by replacing $\alpha$ by $\gamma$.  Then we set $R_g(\alpha, \gamma)=\sum_\rho (\rho, \rho')$, the sum being taken over all relations $\rho$ in which $\alpha$ appears.  In our case,
$(\tilde Q_1 || N)_g$ consists of pairs of the form $(\alpha, \alpha)$ with $\alpha \in \tilde Q_1$.  Therefore $R_g=0$. \\
Now 
$$| \tilde Q_1 ||N| - |(\tilde Q_1 || N)_e| = | (\tilde Q_1 || N)_g \cup (\tilde Q_1 || N)_a |. $$
Hence there exists $\epsilon \geq 0$ such that
$$| \tilde Q_1 ||N| - |(\tilde Q_1 || N)_e| = | (\tilde Q_1 || N)_g | + \epsilon = | \tilde Q_1| + \epsilon .$$
On the other hand, since $(\tilde Q_1 || N)_e$ is a subset of $\tilde Q_1 ||N$, we have that
$$| \tilde Q_1 ||N| - |(\tilde Q_1 || N)_e| = | \tilde Q_1 | + \epsilon \leq | \tilde Q_1 ||N| = |\tilde Q_1| + \epsilon'. $$
Hence $0 \leq \epsilon \leq \epsilon' \leq 2$ and we have 
$$\dim {\HH}^1(B)=1-|\tilde Q_0|+|\tilde Q_1| + \epsilon.$$
So, since $C$ is tilted of type $\tilde{\mathbb{A}}$, its quiver $Q$ contains a unique non-oriented cycle without double zeros and all arrows attached to the cycle either enter it or all leave it, and hence
$$\dim {\HH}^1(B)= 1 - |Q_0|+|Q_1|+|R|+\epsilon=1+|R|+\epsilon.$$ \qed
\end{proof}

\begin{example}
We show that each of the three possible values of $\epsilon \in \{ 0,1,2 \}$ may occur.
\begin{itemize}
\item[(a)] If $C=\tilde C$ is the Kronecker algebra, then ${\HH}^1(\tilde C)=k^3$.  Here $\epsilon = 2$.
\item [(b)] If $C=\tilde C$ is given by the quiver

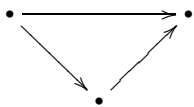
\[ \xymatrix{  \scriptscriptstyle \bullet \ar[rr] \ar[dr] & & \scriptscriptstyle \bullet \\
& \scriptscriptstyle \bullet \ar[ur]
}\] 
then ${\HH}^1(\tilde C)=k^2$.  Here $\epsilon =1$.
\item [(c)] If $B=\tilde C$ is schurian then, since it is monomial, we have 
$${\HH}^1(B)={\HH}^1(C) \oplus k^{|R|} \cong k^{|R|+1}.$$  Here $\epsilon=0$.
\end{itemize}
\end{example}

\noindent
{\bf ACKNOWLEDGEMENTS.} The first author gratefully acknowledges partial
support from the NSERC of Canada and the Universit\'e de Sherbrooke. This paper was started while the second author
was visiting the Universit\'e de Sherbrooke in Qu\'ebec.  She acknowledges
support from the NSERC of Canada, and would like to express
her gratitude to Ibrahim for his warm hospitality.

\end{document}